\documentstyle[12pt]{article}

\font\twlgot =eufm10 scaled \magstep1
\font\egtgot =eufm8
\font\sevgot =eufm7
\font\twlmsb =msbm10 scaled \magstep1
\font\egtmsb =msbm8
\font\sevmsb =msbm7

\newfam\gotfam

\textfont\gotfam\twlgot
\scriptfont\gotfam\egtgot
\scriptscriptfont\gotfam\sevgot

\newfam\msbfam
\textfont\msbfam\twlmsb
\scriptfont\msbfam\egtmsb
\scriptscriptfont\msbfam\sevmsb
\def\Bbb{\protect\pBbb}
\def\pBbb{\relax\ifmmode\expandafter\Bb\else\typeout{You cann't use
Bbb in text mode}\fi}
\def\Bb #1{{\fam\msbfam\relax#1}}

\def\thebibliography#1{\section*{References}\list
   {[\arabic{enumi}]}{\settowidth\labelwidth{#1}\leftmargin\labelwidth
     \advance\leftmargin\labelsep
     \usecounter{enumi}}
     \def\newblock{\hskip .11em plus .33em minus .07em}
     \sloppy\clubpenalty4000\widowpenalty4000
     \sfcode`\.=1000\relax}

\def\op#1{\mathop{\fam0 #1}\limits}

\newcommand{\id}{{\rm Id\,}}

\newcommand{\beq}{\begin{equation}}
\newcommand{\eeq}{\end{equation}}
\newcommand{\ben}{\begin{eqnarray}}
\newcommand{\een}{\end{eqnarray}}
\newcommand{\be}{\begin{eqnarray*}}
\newcommand{\ee}{\end{eqnarray*}}
\newcommand{\bea}{\begin{eqalph}}
\newcommand{\eea}{\end{eqalph}}

\newcommand{\cF}{{\cal F}}

\newcommand{\al}{\alpha}

\newcommand{\bt}{\beta}
\newcommand{\dl}{\delta}
\newcommand{\la}{\lambda}

\newcommand{\f}{\phi}

\newcommand{\Om}{\Omega}

\newcommand{\vt}{\vartheta}
\newcommand{\vf}{\varphi}

\newcommand{\si}{\sigma}

\newcommand{\w}{\wedge}

\newcommand{\dr}{\partial}
\newcommand{\ar}{\op\longrightarrow}

\newcounter{eqalph}
\newcounter{equationa}
\newcounter{theorem}
\newcounter{remark}
\newcounter{proposition}
\newcounter{lemma}
\newcounter{corollary}
\newcounter{definition}

\newenvironment{eqalph}{\stepcounter{equation}
\setcounter{equationa}{\value{equation}}
\setcounter{equation}{0}

\begin{eqnarray}}{\end{eqnarray}\setcounter{equation}{\value{equationa}}}

\def\theremark{\arabic{remark}}

\def\thedefinition{\arabic{definition}}

\newenvironment{proof}{\noindent
{\bf Proof.}}{\hfill $\Box$ \medskip}

\newenvironment{theo}{\refstepcounter{definition}
\bigskip\noindent{\bf Theorem \thedefinition.} \it}{\medskip}
\newenvironment{prop}{\refstepcounter{definition}
\bigskip\noindent{\bf Proposition \thedefinition.}\it}{\medskip}

\newcommand{\mar}[1]{}

\begin{document}
\hbox{}

{\parindent=0pt

{\large\bf The Liouville--Arnold--Nekhoroshev theorem

for non-compact
invariant manifolds}
\bigskip

{\bf Emanuele Fiorani}$^1$, {\bf Giovanni Giachetta}$^2$ {\bf and 
Gennadi Sardanashvily}$^3$
\bigskip

\begin{small}

$^1$ Department of Mathematics "F.Enriques", University of Milano,
20133 Milano, Italy

E-mail: fiorani@mat.unimi.it

$^2$ Department of Mathematics and Informatics, University of Camerino, 62032
Camerino (MC), Italy

E-mail: giovanni.giachetta@unicam.it

$^3$ Department of Theoretical Physics, Physics Faculty, Moscow State
University, 117234 Moscow, Russia

E-mail: sard@grav.phys.msu.su, URL: 
http//webcenter.ru/$\sim$sardan
\medskip

\end{small}


\bigskip

{\bf Abstract.}

Under certain conditions, generalized action-angle coordinates can be
introduced near non-compact invariant manifolds of completely and
partially integrable Hamiltonian systems.
\medskip

PACS numbers: 45.20.Jj, 02.30.Ik


}

\section{Introduction}

Let us recall that an autonomous  Hamiltonian system on a $2n$-dimensional
symplectic manifold
is said to be completely integrable  if
there exist $n$ independent integrals of motion in involution. By
virtue of the
classical Liouville--Arnold theorem \cite{arn,laz}, such a  system admits
action-angle coordinates around a connected regular compact invariant manifold.
In a more general setting, one considers Hamiltonian systems  having
partial integrability, i.e., $k\leq n$ independent integrals of 
motion in involution.
The Nekhoroshev theorem for these systems
\cite{gaeta,nekh}  generalizes  both the Poincar\'e--Lyapunov
theorem ($k=1$) and the above mentioned
Liouville--Arnold theorem ($k=n$). The Nekhoroshev theorem
in fact falls into
two parts. The first part states the sufficient conditions for
an open neighbourhood of
an invariant torus $T^k$ to be a trivial fibre bundle (see \cite{gaeta}
for a detailed exposition). The second one provides this bundle with 
partial action-angle
coordinates similarly to the case of complete
integrability.

The present work addresses completely and partially integrable 
Hamiltonian systems
whose invariant manifolds need not be compact. This is the case of
any autonomous Hamiltonian
system because its Hamiltonian, by definition, is an integral of motion.
In the preceding papers, we have shown that, if an open neighbourhood
of a non-compact invariant
manifold of a completely integrable
Hamiltonian system is a trivial bundle, it can
be equipped with the generalized action-angle coordinates which bring
a symplectic form into the canonical form \cite{acang1,acang2}.
Here, we prove that, under certain conditions,
an open neighbourhood of a regular non-compact invariant
manifold of a completely integrable system is a trivial bundle (see 
parts (A) -- (C) in
the proof of Theorem \ref{dd1}) and, consequently,
it can be equipped with the generalized action-angle
coordinates (see part (D) of this theorem).
Then, this result is extended to partially integrable Hamiltonian systems
(see Theorem \ref{dd3}). The proof of Theorem \ref{dd3} mainly follows that
of Theorem \ref{dd1}. Note that part (D) in
the proof of
Theorem \ref{dd1} can be simplified by
the choice of a Lagrangian
section $\si$, but this is not the case of partially integrable systems.
This proof also shows that, from the beginning, one can separate integrals of
motion whose trajectories live in tori.

It should be emphasized that the results of Theorem \ref{dd3} are not
limited by the scope of autonomous mechanics. Any time-dependent Hamiltonian
system of $n$ degrees of
freedom can be
extended to an autonomous Hamiltonian system of $n+1$ degrees of freedom
which has at least one integral of motion, namely, its Hamiltonian 
\cite{acang1}.
Thus, any time-dependent Hamiltonian
system can be seen as a partially integrable autonomous Hamiltonian 
system whose invariant
manifolds are never compact because of the time axis.
Just the time is a generalized angle coordinate corresponding
to a Hamiltonian of this autonomous system.

One also finds reasons in quantum theory in order to introduce
generalized action-angle variables. In particular, quantization with respect
to these variables enables one to include a Hamiltonian in
the quantum algebra \cite{acang1,quant}.

\section{Completely integrable systems}

Let $(Z,\Om)$ be a $2n$-dimensional symplectic manifold, and let it
admit $n$ real smooth functions $\{F_\la\}$, which are pairwise in involution
and independent almost everywhere on $Z$. The latter implies that
the set of non-regular points, where
the morphism
\mar{d1}\beq
\pi=\op\times^\la F_\la : Z\to\Bbb R^n \label{d1}
\eeq
fails to be a submersion, is nowhere dense.
Bearing in mind physical applications,
we agree to think of one of the functions $F_\la$ as being a Hamiltonian
and of the other as first integrals of motion. Accordingly,
their common level surfaces are called
invariant surfaces.

Let $M$ be a regular invariant surface, i.e.,
the morphism $\pi$ (\ref{d1}) is a submersion at all
points of $M$ or, equivalently, the $n$-form $\op\w^\la dF_\la$ 
vanishes nowhere on $M$.
Hence, $M$ is a closed imbedded submanifold of $Z$. There exists its open
neighbourhood $U$ such
that the morphism $\pi$ is a submersion on $U$, i.e.,
\mar{d3}\beq
\pi: U\to N=\pi(U) \label{d3}
\eeq is a
fibred manifold over an open subset $N\subset\Bbb R^n$.
The vertical tangent bundle $VU$ of
$U\to N$ coincides with the $n$-dimensional distribution on $U$
spanned by the
Hamiltonian vector fields $\vt_\la$ of the functions $F_\la$.
Integral manifolds of this distribution are components
of the fibres of $\pi$. They are Lagrangian submanifolds of $Z$.
Let $U$ be connected. Then $N$ is a domain. Without loss of
generality, one can suppose that
there exists a section of $U\to N$.

If $M$ is connected and compact,
we come to the conditions of the Liouville--Arnold theorem.
If $M$ need not be compact, one should require something more.

\begin{theo} \mar{dd1} \label{dd1}
Let $M$ be a connected regular invariant manifold of a completely
integrable Hamiltonian system $\{F_\la\}$ and let $U$ be an open 
neighbourhood as
above. Let us additionally assume that: (i)  all fibres of the fibred
manifold $U\to N$ (\ref{d3}) are mutually diffeomorphic, (ii)
the Hamiltonian vector fields $\vt_\la$ on $U$ are complete.
Then, there exists a domain $N$ so that $U\to N$ is a trivial bundle
\mar{z10}\beq
U=N\times(\Bbb R^{n-m}\times T^m), \label{z10}
\eeq
provided with the generalized
action-angle coordinates
$(I_\la;x^a; \f^i)$ such that the integrals of motion $F_\la$ depend 
only on the action
coordinates $I_\al$ and the symplectic form $\Om$ on $U$ reads
\mar{d10}\beq
\Om=dI_a\w dx^a +dI_i\w d\f^i. \label{d10}
\eeq
\end{theo}

\begin{proof} (A)
Since Hamiltonian vector fields $\vt_\al$ on $U$
are complete and mutually commutative, their flows assemble into the 
additive Lie group
$\Bbb R^n$. This group is naturally identified with its Lie
algebra, and its group space is a
vector space coordinated by parameters $(s^\la)$
of the flows with respect to the basis $\{e_\la\}$ for its Lie algebra.
This group acts in $U$  so
that its generators $e_\la$ are represented by
the Hamiltonian vector fields $\vt_\la$ and
its orbits are fibres of the fibred manifold $U\to N$.
Given a point $r\in N$, the action of $\Bbb R^n$
in the fibre $M_r=\pi^{-1}(r)$ factorizes as
\mar{d4}\beq
\Bbb R^n\times M_r\to G_r\times M_r\to M_r \label{d4}
\eeq
through the free transitive
action in $M_r$ of the factor group $G_r=\Bbb R^n/K_r$, where $K_r$ 
is the isotropy group of
an arbitrary point of $M_r$. It is the same group for all points
because $\Bbb R^n$ is an Abelian group.
Since the fibres $M_r$ are mutually diffeomorphic, all isotropy groups $K_r$
are isomorphic to the group $\Bbb Z^m$ for some fixed $m$, $0\leq m\leq n$,
and the groups
$G_r$ are isomorphic to the additive group $\Bbb R^{n-m}\times T^m$.
Let us show that the fibred manifold $U\to N$ (\ref{d3}) is a 
principal bundle with
the structure group $G_0$, where we denote $\{0\}=\pi(M)$. For this purpose,
let us determine
isomorphisms $\rho_r: G_0\to G_r$ of the group $G_0$ to the groups 
$G_r$, $r\in N$.
Then, a desired fibrewise action of $G_0$ in $U$ is given by the law
\mar{d5}\beq
G_0\times M_r\to\rho_r(G_0)\times M_r\to M_r. \label{d5}
\eeq

(B) Generators of each isotropy subgroup
$K_r$ of $\Bbb R^n$ are given by $m$ linearly independent vectors of the group
space $\Bbb R^n$. One can show that there are ordered collections of generators
  $(v_1(r),\ldots,v_m(r))$
of the groups $K_r$ such that $r\mapsto v_i(r)$
are smooth $\Bbb R^n$-valued fields on $N$. Indeed, given a vector $v_i(0)$
and a section $\si$ of the fibred manifold $U\to N$,
each field $v_i(r)=(s^\al(r))$ is the unique smooth solution of the equation
\mar{d40}\beq
g(s^\al)\si(r)=\si(r), \qquad  (s^\al(0))=v_i(0), \label{d40}
\eeq
on an open neighbourhood of $\{0\}$. Without loss of generality,
one can assume that this neighbourhood is $N$.
Let us consider the decomposition
\be
v_i(0)=B_i^a(0) e_a + C_i^k(0) e_k, \qquad a=1,\ldots,n-m, \qquad 
k=1,\ldots, m,
\ee
where $C_i^k(0)$ is a non-degenerate matrix.
Since the fields $v_i(r)$ are smooth, there exists an open 
neighbourhood of $\{0\}$,
say $N$ again, where the matrices $C_i^k(r)$ remain non-degenerate. 
Then, there is
a unique linear
morphism
\mar{d6}\beq
A_r=\left(
\begin{array}{ccc}
\id & \qquad & (B(r)-B(0))C^{-1}(0) \\
0 & & C(r)C^{-1}(0)
\end{array}
\right) \label{d6}
\eeq
of the vector space $\Bbb R^n$ which transforms
its frame $v_\al(0)=\{e_a,v_i(0)\}$
into the frame $v_\al(r)=\{e_a,v_i(r)\}$.
Since it is also an automorphism of the group $\Bbb R^n$
sending $K_0$ onto $K_r$, we obtain a desired isomorphism
$\rho_r$ of the group $G_0$ to the group $G_r$. Let an element $g$ of 
the group $G_0$
be the coset of an element $g(s^\la)$ of the group $\Bbb R^n$. Then, it
acts in $M_r$ by the rule (\ref{d5}) just as the element 
$g((A_r^{-1})^\la_\bt s^\bt)$
of the group $\Bbb R^n$ does. Since entries of the matrix $A$ (\ref{d6}) are
smooth functions on $N$, this action of the group $G_0$ in $U$ is 
smooth. It is free, and
$U/G_0=N$. Then, the fibred manifold $U\to N$ is a principal bundle 
with the structure
group $G_0$ which is trivial because $N$ is a domain.

(C) Given a section $\si$ of the principal bundle $U\to N$, its 
trivialization $U=N\times G_0$
is defined by assigning the points $\rho^{-1}(g_r)$ of the
group space $G_0$ to points
$g_r\si(r)$, $g_r\in G_r$, of a fibre $M_r$.
Let us endow $G_0$ with the standard
coordinate atlas
$(y^\la)=(t^a;\vf^i)$ of the group $\Bbb R^{n-m}\times T^m$. We 
provide $U$ with a desired
trivialization (\ref{z10}) with respect to the coordinates $(J_\la;t^a;\vf^i)$,
where $J_\la(u)=F_\la(u)$, $u\in U$, are coordinates on the base $N$.
The Hamiltonian vector fields $\vt_\la$ on $U$ relative to these 
coordinates read
\mar{ww25}\beq
\vt_a=\dr_a, \qquad \vt_i=-(BC^{-1})^a_i\dr_a +
(C^{-1})_i^k\dr_k.\label{ww25}
\eeq
In particular, the Hamilton equation takes the form
\be
\dot J_\la=0, \qquad \dot y^\la= f^\la(J_\al).
\ee

(D) Since fibres of $U\to N$ are Lagrangian manifolds,
the symplectic form $\Om$ on $U$ is given by the coordinate expression
\mar{ac1}\beq
\Om=\Om^{\al\bt}dJ_\al\w dJ_\bt + \Om^\al_\bt dJ_\al\w dy^\bt. \label{ac1}
\eeq
Let us bring it into the canonical form (\ref{d10}).
The Hamiltonian vector fields $\vt_\la$ obey the relations 
$\vt_\la\rfloor\Om=-dJ_\la$,
which take the coordinate form
\mar{ww22}\beq
\Om^\al_\bt \vt^\bt_\la=\dl^\al_\la. \label{ww22}
\eeq
It follows that $\Om^\al_\bt$ is a nondegenerate matrix whose entries
are independent of coordinates $y^\la$.
By virtue of
the well-known K\"unneth formula for the de Rham cohomology of 
manifold product,
the closed form $\Om$ (\ref{ac1}) on $U$ (\ref{z10})
is exact, i.e., $\Om=d\Xi$ where $\Xi$ reads
\mar{ac2}\beq
\Xi=\Xi^\al(J_\la,y^\la)dJ_\al + \Xi_i(J_\la) d\vf^i. \label{ac2}
\eeq
Because entries of $d\Xi=\Om$ are independent of $y^\la$, we obtain 
the following.

(i) $\Om^\la_i=\dr^\la\Xi_i-\dr_i\Xi^\la$. Consequently, $\dr_i\Xi^\la$ are
independent of $\vf^i$, i.e., $\Xi^\la$ are at most affine in $\vf^i$ 
and, therefore,
are independent of $\vf^i$ since these are cyclic coordinates. Hence,
$\Om^\la_i=\dr^\la\Xi_i$ and $\dr_i\rfloor\Om=-d\Xi_i$. A glance at 
the last equality
shows that $\dr_i$ are Hamiltonian vector fields. It follows that
we can substitute $m$ integrals of motion  among $F_\la$ withe the functions
$\Xi_i$, which we continue to denote $F_i$. The Hamiltonian vector 
fields of these new $F_i$
are tangent to invariant tori. In this case, the matrix $B$ in the expressions
(\ref{d6}) and (\ref{ww25}) is the zero one, and the Hamiltonian vector fields
$\vt_\la$ read
\mar{ww25'}\beq
\vt_a=\dr_a, \qquad \vt_i=(C^{-1})_i^k\dr_k. \label{ww25'}
\eeq
Moreover, the coordinates $t^a$ are exactly the flow parameters $s^a$.
Substituting the expressions (\ref{ww25'}) into the conditions 
(\ref{ww22}), we obtain
\be
\Om=\Om^{\al\bt}dJ_\al\w dJ_\bt +dJ_a\w ds^a + C^i_k dJ_i\w d\vf^k.
\ee
It follows that $\Xi_i$ are independent of $J_a$, and so are 
$C^k_i=\dr^k\Xi_i$.

(ii) $\Om^\la_a=-\dr_a\Xi^\la=\dl^\la_a$. Hence, $\Xi^a=-s^a+E^a(J_\la)$
and $\Xi^i$ are independent of $s^a$.

In view of items (i) -- (ii), the Liouville form $\Xi$ (\ref{ac2}) reads
\be
\Xi=(-s^a+E^a(J_\la))dJ_a + E^i(J_\la)dJ_i + \Xi_i(J_j) d\vf^i.
\ee
Since the matrix $\dr^k\Xi_i$ is nondegenerate,
we perform the coordinate transformation $I_a=J_a$, $I_i=\Xi_i(J_j)$, 
and obtain
\be
\Xi=(-s^a+E'^a(I_\la))dI_a + E'^i(I_\la) dI_i + I_i d\vf^i.
\ee
Finally, put
\mar{ee2}\beq
x^a=s^a-E'^a(I_\la). \qquad
\f^i=\vf^i-E'^i(I_\la) \label{ee2}
\eeq
in order to obtain the desired action-angle coordinates
\be
I_a=J_a, \qquad I_i(J_j), \qquad x^a(J_\la,s^a), \qquad \f^i(J_\la,\vf^k).
\ee
The shifts (\ref{ee2}) correspond to the choice of a Lagrangian section $\si$.
\end{proof}

Let us remark that the generalized action-angle coordinates in 
Theorem \ref{dd1} are by no
means unique. For instance, the canonical coordinate transformations
\mar{ww26}\beq
I_a=f_a(I'_\la), \qquad I_i=I'_i, \qquad x'^a=\frac{\dr f_b}{\dr 
I'_a}x^b, \qquad \f'^i= \f^i +
\frac{\dr f_a}{\dr I'_i}x^a. \label{ww26}
\eeq
give new generalized action-angle coordinates on $U$.

\section{Partially integrable systems}

Let a $2n$-dimensional symplectic manifold $(Z,\Om)$ admit $k< n$
smooth real functions $\{F_\la\}$,
which are pairwise in involution
and independent almost everywhere on $Z$. Let us consider the morphism
\mar{d11}\beq
\pi=\op\times^\la F_\la : Z\to\Bbb R^k, \label{d11}
\eeq
and its regular connected common level surface $W$. There exists an 
open connected
neighbourhood $U_W$ of $W$ such that
\mar{d13}\beq
\pi: U_W\to V_W=\pi(U_W) \label{d13}
\eeq
is a fibred manifold over a domain $V_W$ in $\Bbb R^k$.
Restricted to $U_W$, the Hamiltonian vector fields $\vt_\la$
of functions $F_\la$ define a $k$-dimensional distribution
and the corresponding regular foliation $\cF$ of $U_W$.
Its leaves are isotropic. They are
located in fibres of the fibred manifold $U_W\to V_W$
and, moreover, make up regular foliations of these fibres.

Let us assume that the foliation $\cF$ has a total transversal
manifold $S$ and its holonomy pseudogroup on $S$ is trivial. Then,
$\cF$ is a fibred manifold
\mar{d20}\beq
\pi_1: U_W\to S' \label{d20}
\eeq
and $S=\si(S')$ is its section \cite{mol}.
Thereby, the fibration $\pi$ (\ref{d13}) factorizes as
\be
\pi: U_W\ar^{\pi_1} S'\ar^{\pi_2} V_W
\ee
through
the fibration $\pi_1$ (\ref{d20}). The map $\pi_2$ reads 
$\pi_2=\pi\circ \si$ and,
consequently, it is also a fibred manifold.

\begin{prop} \mar{dd2} \label{dd2}
Let us assume that there  exists a domain $N\subset S'$ such that:
(i) the fibres of the fibred manifold $\pi_1$ (\ref{d20}) over $N$ are
mutually diffeomorphic, (ii)
the Hamiltonian vector fields $\vt_\la$ on $U=\pi_1^{-1}(N)$ are complete.
Then, there exists a domain in $S'$, say $N$ again, such that $U\to N$ is a
trivial principal bundle with the structure group
$\Bbb R^{k-m}\times T^m$.
\end{prop}

\begin{proof}
The proof is a straightforward repetition of parts (A) -- (B) in the
proof of Theorem \ref{dd1}.
\end{proof}

Furthermore, one can always choose the
domain $N$ in Proposition \ref{dd2} as the domain of a fibred chart of
$\pi_2$.  Following part (C) in the proof of Theorem \ref{dd1}, we can provide
$U\to N$ with the trivialization
\mar{d21}\beq
U=N\times\Bbb R^{k-m}\times T^m, \label{d21}
\eeq
coordinated by $(J_\la;z^A;y^\la)$ where:
(i) $J_\la(u)=F_\la(u)$, $u\in U$, are coordinates on the base $V$, (ii)
$(J_\la;z^A)$ are coordinates on $N$, and (iii)
$(y^\la)=(t^a;\vf^i)$ are coordinates on $\Bbb R^{k-m}\times T^m$.
The Hamiltonian vector fields $\vt_\al$ on $U$ with respect to these 
coordinates read
\mar{d22}\beq
\vt_a=\dr_a, \qquad  \qquad  \vt_i=\vt_i^a(J_\la,z^A)\dr_a 
+\vt_i^k(J_\la,z^A)\dr_k.\label{d22}
\eeq
Since fibres of $U\to N$ are isotropic,
the symplectic form $\Om$ on $U$
relative to the
coordinates $(J_\la;z^A;y^\la)$ is given by the expression
\mar{d23}\beq
\Om=\Om^{\al\bt}dJ_\al\w dJ_\bt + \Om^\al_\bt dJ_\al\w dy^\bt +
\Om_{AB}dz^A\w dz^B +\Om_A^\la dJ_\la\w dz^A +
  \Om_{A\bt} dz^A\w dy^\bt. \label{d23}
\eeq
The Hamiltonian vector fields $\vt_\la$ obey the relations 
$\vt_\la\rfloor\Om=-dJ_\la$,
which give the conditions
\be
\Om^\al_\bt \vt^\bt_\la=\dl^\al_\la, \qquad \Om_{A\bt}\vt^\bt_\la=0.
\ee
The first of them shows that $\Om^\al_\bt$ is a non-degenerate matrix
independent of coordinates $y^\la$. Then, the second one implies 
$\Om_{A\bt}=0$.
The rest is a minor modification of part (D) in the proof of Theorem \ref{dd1}.

The symplectic form $\Om$ (\ref{d23})
on $U$ is exact, and the Liouville form is
\be
\Xi=\Xi^\al(J_\la,z^B,y^\la)dJ_\al + \Xi_i(J_\la,z^B) d\vf^i
+\Xi_A(J_\la,z^B,y^\la)dz^A.
\ee
Since $\Xi_a=0$ and $\Xi_i$ are independent of $\vf^i$, one easily obtains from
the relations $\Om_{A\bt}=\dr_A\Xi_\bt-\dr_\bt\Xi_A=0$ that
$\Xi_i$ are  independent of coordinates $z^A$,
while $\Xi_A$ are independent of coordinates $y^\la$. Hence,
the Liouville form reads
\be
\Xi=\Xi^\al(J_\la,z^B,y^\la)dJ_\al + \Xi_i(J_\la) d\vf^i
+\Xi_A(J_\la,z^B)dz^A
\ee
(cf. (\ref{ac2})). Running through item (i), we observe that, in the 
case of a partially integrable
system, one can also separate integrals of motion $F_i$ whose 
Hamiltonian vector
fields are tangent to invariant tori. Then, the Hamiltonian vector 
fields (\ref{d22})
take the form
\be
\vt_a=\dr_a, \qquad  \qquad  \vt_i=\vt_i^k(J_\la,z^A)\dr_k.
\ee
Following items (i) -- (ii) of part (D), we obtain
\be
\Xi=(-s^a+E^a(J_\la,z^B))dJ_a + E^i(J_\la,z^B) dJ_i + \Xi_i(J_j) d\vf^i +
\Xi_A(J_\la,z^B)dz^A.
\ee
Finally, the coordinates
\be
x^a=-s^a+E^a(J_\la,z^B), \qquad I_i=\Xi_i(J_j), \qquad I_a=J_a,
\qquad \f^i=\vf^i-E^j(J_\la,z^B)\frac{\dr J_j}{\dr I_i}
\ee
bring $\Om$ into the form
\mar{d26}\beq
\Om= dI_a\w dx^a + dI_i\w d\f^i +\Om_{AB}(I_\la,z^B) dz^A\w dz^B 
+\Om_A^\la(I_\la,z^B) dI_\la\w dz^A.
\label{d26}
\eeq
Therefore, one can think of these coordinates as being partial 
generalized action-angle coordinates.
The Hamiltonian vector fields of integral of motions with respect to 
these coordinates read
\be
\vt_a=\frac{\dr}{\dr x^a}, \qquad \vt_i=\frac{\dr J_i}{\dr 
I_j}\frac{\dr}{\dr \f^j}.
\ee

Thus, we have proved the following.

\begin{theo}  \label{dd3} \mar{dd3}
Given a partially integrable Hamiltonian system $\{F_\la\}$ on a 
symplectic manifold $(Z,\Om)$,
let $W$ be its regular connected level surface, and let $M\subset W$ 
be a leaf of the
characteristic foliation $\cF$ of the distribution
generated by the Hamiltonian vector fields $\vt_\la$ of
$F_\la$. Let $M$ have an open satured neighbourhood $U\subset Z$ such that:
(i) the foliation $\cF$ of $U$ admits a transversal manifold $S$
and its holonomy pseudogroup on $S$ is trivial, (ii) the leaves of 
this foliation
are
mutually diffeomorphic, (iii) Hamiltonian vector fields $\vt_\la$
on $U$ are complete. Then, there exists an open satured neighbourhood of $M$,
say $U$ again, which is a trivial
bundle (\ref{d21}), provided with the particular
coordinates $(I_\la;z^A;x^a;\f^i)$ such that
the integrals of motion $F_\la$ depend only on the
coordinates $I_\al$ and the symplectic form $\Om$ on $U$
is brought into the form (\ref{d26}).
\end{theo}


\begin{thebibliography}{}


\bibitem{arn} Arnold V (Ed.) 1988 {\it Dynamical Systems III}
(Berlin: Springer-Verlag)

\bibitem{acang1} Fiorani F, Giachetta G and Sardanashvily G 2002 {\it 
J. Math. Phys.}
{\bf 43} 5013

\bibitem{gaeta} Gaeta G 2002 {\it Ann. Phys.} {\bf 297} 157

\bibitem{acang2} Giachetta G, Mangiarotti L and Sardanashvily G 2002
{\it J. Phys. A} {\bf 35} L439

\bibitem{quant} Giachetta G, Mangiarotti L and Sardanashvily G 2002
{\it Phys. Lett. A} {\bf 301} 53

\bibitem{laz} Lazutkin V 1993 {\it  KAM Theory and Semiclassical
Approximations to  Eigenfunctions} (Berlin: Springer-Verlag)

\bibitem{mol} Molino P 1988 {\it Riemannian Foliations} (Boston:
Birkh\"auser)

\bibitem{nekh} Nekhoroshev N 1994
{\it Funct. Anal. Appl.} {\bf 28} 128

\end{thebibliography}
\end{document}